\documentclass{amsart}

\newtheorem{theorem}{Theorem}[section]

\newtheorem{lemma}[theorem]{Lemma}
\newtheorem{corollary}[theorem]{Corollary}

\newtheorem{example}[theorem]{Example}

\newcommand{\dH}{\Delta_\mathcal{H}}
\newcommand{\dAH}{\Delta_{\mathcal{A},\mathcal{H}}}
\newcommand{\charac}{\chi}
\newcommand{\A}{\mathcal{A}}
\newcommand{\h}{\mathcal{H}}
\newcommand{\Sn}{\mathcal{S}_n}
\newcommand{\Bn}{\mathcal{B}_n}
\newcommand{\N}{\mathbb{N}}
\newcommand{\Z}{\mathbb{Z}}
\newcommand{\R}{\mathbb{R}}
\newcommand{\LH}{L_\mathcal{H}}
\newcommand{\G}{\widehat{G}}
\newcommand{\Kn}{\widehat{K_n}}
\newcommand{\dG}{\Delta_{\G,\Sn}}
\newcommand{\hrev}{\overline{h}}
\newcommand{\f}{\mathcal{F}}
\newcommand{\ph}{P_\h}

\newcommand{\hilb}{\text{Hilb}}
\newcommand{\J}{\mathcal{J}}

\newcommand{\tail}{T}
\newcommand{\Hn}{\hat{0}}
\newcommand{\chrom}{P}
\begin{document}
\title{Link complexes of subspace arrangements}
\thanks{This research was carried out, and most of the article was written, while the author was a postdoc at Institut Mittag-Leffler.}
\author{Axel Hultman}
\address{Department of Mathematics, KTH, SE-100 44 Stockholm, Sweden }
\email{axel@math.kth.se}

\begin{abstract}
Given a simplicial hyperplane arrangement $\h$ and a subspace
arrangement $\A$ embedded in $\h$, we define a simplicial complex
$\dAH$ as the subdivision of the link of $\A$ induced by $\h$. In particular, this generalizes Steingr\'{\i}msson's coloring complex of a graph.

We do the following:
\begin{enumerate}
\item When $\A$ is a hyperplane arrangement, $\dAH$ is shown to be shellable. As a special case, we answer affirmatively a question of Steingr\'{\i}msson on coloring complexes.
\item For $\h$ being a Coxeter arrangement of type $A$ or $B$ we obtain a close connection between the Hilbert series of the Stanley-Reisner ring of $\dAH$ and the characteristic polynomial of $\A$. This extends results of Steingr\'{\i}msson and provides an interpretation of chromatic polynomials of hypergraphs and signed graphs in terms of Hilbert polynomials.
\end{enumerate}
\end{abstract}

\maketitle

\section{Introduction}
In \cite{steingrimsson}, Steingr\'{\i}msson introduced the {\em coloring complex} $\Delta_G$. This is a simplicial complex associated with a graph $G$. The Hilbert polynomial of its Stanley-Reisner ring $k[\Delta_G]$ is closely related to the chromatic polynomial $\chrom_G(x)$ in a way that is made precise in Section \ref{se:h-vector}. 

Answering a question of Steingr\'{\i}msson, Jonsson \cite{jonsson} proved that $\Delta_G$ is a Cohen-Macaulay complex by showing that it is constructible. In particular, $\Delta_G$ being Cohen-Macaulay imposes restrictions on the Hilbert polynomial of $k[\Delta_G]$, hence on $\chrom_G(x)$. 

Since $\Delta_G$ is a Cohen-Macaulay complex, a natural question, asked already in \cite{steingrimsson}, is whether it is shellable --- a stronger property than constructibility.

In \cite{steingrimsson}, $\Delta_G$ was defined in a combinatorially very explicit way. Another way to view $\Delta_G$ is, however, as a simplicial decomposition of the link (i.e.\ intersection with the unit sphere) of the {\em graphical hyperplane arrangement} associated with $G$. In this guise, $\Delta_G$ appeared in work of Herzog, Reiner and Welker \cite{HRW}. Adopting this point of view, one may define a similar complex $\dAH$ for any subspace arrangement $\A$, as long as it has an embedding in a simplicial hyperplane arrangement $\h$.

This paper has two goals. The first is addressed in Section \ref{se:shellable} where we show that $\dAH$ is shellable whenever $\A$ consists of hyperplanes. In particular, this proves that the coloring complexes are shellable. 

The chromatic polynomial of $G$ is essentially the characteristic polynomial of the corresponding graphical hyperplane arrangement. Bearing this in mind, one may hope to extend the aforementioned connection between the Hilbert polynomial of $k[\Delta_G]$ and $\chrom_G(x)$ to more general complexes $\dAH$. Achieved in Section \ref{se:h-vector}, our second goal is to carry out this extension whenever $\h$ is a Coxeter arrangement of type $A$ or $B$. When $\A$ consists of hyperplanes and $\h$ is of type $A$, Steingr\'{\i}msson's result is recovered.

We define the complexes $\dAH$ in Section \ref{se:defin} after reviewing some necessary background in the next section.

\section{Preliminaries}
\subsection{Subspace arrangements and characteristic polynomials}
By the term {\em subspace arrangement} we mean a finite collection
$\A=\{A_1, \dots, A_t\}$ of linear subspaces, none of which contains
another, of some ambient vector space. In our case, the ambient space
will always be $\R^n$ for some $n$. To $\A$ we associate the {\em
  intersection lattice} $L_\A$ which consists of all intersections of
subspaces in $\A$ ordered by reverse inclusion. (We emphasize the fact
that $\A$ contains no strictly affine subspaces; in particular
this implies that $L_\A$ is indeed a lattice.)

An important invariant of the arrangement $\A$ is its {\em characteristic polynomial} 
\[
\charac(\A;x) = \sum_{Y\in L_\A}\mu(\Hn,Y)x^{\dim (Y)},
\]
where $\mu$ is the M\"obius function of $L_\A$ and $\Hn = \R^n$ is the smallest element in $L_\A$.

Given a subspace $A \in \A$, we define two new arrangements, namely the {\em deletion}
\[
\A\setminus A = \A \setminus \{A\}
\]
and the {\em restriction}
\[
\A/A = \max \{A\cap B| B\in \A\setminus A\},
\]
where $\max \mathcal{S}$ denotes the collection of inclusion-maximal
members of a set family $\mathcal{S}$. Another way to think of $\A/A$
is as the set of elements covering $A$ in $L_\A$. In this way, we may
extend the definition of $\A/A$ to arbitrary $A\in L_\A$. We consider $\A\setminus A$ to be an arrangement in $\R^n$, whereas $\A/A$ is an arrangement in $A$.

When $\A$ is a hyperplane arrangement, the next result is
standard. We expect the general case to be known, too, although we have been unable to find it in the literature.
\begin{theorem}[Deletion-Restriction]\label{th:DR}
For a subspace arrangement $\A$ and any subspace $A\in \A$, we have
\[
\charac(\A;x) = \charac(\A\setminus A;x)-\charac(\A/A;x).
\]
\begin{proof}
Choose $Y \in L_\A$. We claim that
\[
\mu_\A(\Hn,Y)=\begin{cases}
\mu_{\A\setminus A}(\Hn,Y)-\mu_\A(A,Y) & \text{ if $Y\in L_{\A\setminus A}$,}\\
-\mu_\A(A,Y) & \text{ otherwise,}
\end{cases}
\]
where $\mu_\A$ denotes the M\"obius function of $L_\A$ which we think of as a function $L_\A \times L_\A \to \Z$ with $S\not \leq T \Rightarrow \mu_\A(S,T)=0$ (and similarly for $\A\setminus A$). 

The claim is true if $Y = \Hn = \R^n$, so assume it has been verified for all $Z< Y$ in $L_\A$. If $Y \in L_{\A\setminus A}$ we obtain
\[
\begin{split}
\mu_\A(\Hn,Y) & = - \sum_{\Hn\leq Z < Y} \mu_\A(\Hn,Z) =
  -\sum_{\buildrel {\Hn
  \leq Z < Y} \over {Z \in L_{\A\setminus A}}} \mu_{\A\setminus
  A}(\Hn,Z)+\sum_{A \leq Z < Y}\mu_\A(A,Z) \\
& = \mu_{\A\setminus A}(\Hn,Y) - \mu_\A(A,Y),
\end{split}
\]
as desired. If, on the other hand, $Y\not \in L_{\A\setminus A}$, then
there is a unique largest element in $L_{\A\setminus A}$ which is
below $Y$ in $L_\A$, namely the join of all atoms (weakly) below $Y$ except
$A$; call this element $W$. If $W=\Hn$, then $Y=A$ and we are done. Otherwise,
\[
\begin{split}
\mu_\A(\Hn,Y) & = - \sum_{\Hn\leq Z < Y} \mu_\A(\Hn,Z) =
  -\sum_{\buildrel {\Hn
  \leq Z \leq W} \over {Z \in L_{\A\setminus A}}} \mu_{\A\setminus
  A}(\Hn,Z)+\sum_{A \leq Z < Y}\mu_\A(A,Z) \\
& =  \sum_{A \leq Z < Y}\mu_\A(A,Z) = - \mu_\A(A,Y),
\end{split}
\]
establishing the claim.

We conclude that
\[
\charac(\A;x)=\sum_{Y\in L_{\A\setminus A}}\mu_{\A\setminus
    A}(\Hn,Y)x^{\dim(Y)} - \sum_{Y\geq A}\mu_{\A}(A,Y)x^{\dim(Y)}.
\]
Not every $Y$ in the last sum belongs to
$L_{\A/A}$ in general; the latter is join-generated by the elements
covering $A$ in $L_\A$. However, it follows from Rota's Crosscut
theorem \cite{rota} that for every $Y\geq A$ in $L_A$,
\[
\mu_\A(A,Y) = \begin{cases}
\mu_{\A/A}(A,Y) & \text{ if $Y\in L_{\A/A}$,}\\
0 & \text{ otherwise.}
\end{cases}
\]
Thus, 
\[
\sum_{Y\geq A}\mu_{\A}(A,Y)x^{\dim(Y)} = \charac(\A/A;x),
\]
and the theorem follows.
\end{proof}
\end{theorem}
Two (families of) hyperplane arrangements are of particular importance
to us. The first is the {\em braid arrangement} $\Sn$. This is an
arrangement whose ambient space is $\{(x_1, \dots, x_n)\in \R^n\mid x_1+\dots + x_n
  = 0\} \cong \R^{n-1}$. The $\binom{n}{2}$ hyperplanes in $\Sn$ are given by the equations $x_i = x_j$ for all $1\leq i<j\leq n$.

The braid arrangement is the set of reflecting hyperplanes of a Weyl
group of type $A$. Considering type $B$ instead, we find our second
important family of arrangements. Explicitly, $\Bn$ is the arrangement
of the $n^2$ hyperplanes in $\R^n$ that are given by the equations
$x_i = \tau x_j$ for all $1\leq i < j \leq n$, $\tau \in \{-1,1\}$,
and $x_i=0$ for all $i\in [n]=\{1, \dots, n\}$.   

\subsection{Stanley-Reisner rings and $h$-polynomials}
Let $\Delta$ be a simplicial complex on the vertex set $[n]$. Regarding the vertices as variables, we want to consider the ring of polynomials that live on $\Delta$. To this end, for a field $k$, we define the {\em Stanley-Reisner ideal} $I_\Delta \subseteq k[x_1, \dots, x_n]$ by
\[
I_\Delta = \langle \{x_{i_1}\dots x_{i_t}| \{i_1, \dots, i_t\} \not \in \Delta\}\rangle. 
\]
The quotient ring
\[
k[\Delta] = k[x_1, \dots, x_n]/I_\Delta
\]
is the {\em Stanley-Reisner ring} of $\Delta$, which is a graded algebra with the standard grading by degree. When speaking of algebraic properties, such as Cohen-Macaulayness, of $\Delta$ we have the corresponding properties of $k[\Delta]$ in mind. 

Given a simplicial complex $\Delta$ of dimension $d-1$, its {\em $h$-polynomial} is 
\[
h(\Delta;x) = \sum_{i=0}^d f_{i-1}(x-1)^{d-i},
\]
where $f_i$ is the number of $i$-dimensional simplices in $\Delta$ (including $f_{-1} = 1$ if $\Delta$ is nonempty). One important feature of the $h$-polynomial is that it carries all information needed to compute the Hilbert series of $k[\Delta]$. Specifically,
\[
\hilb(k[\Delta];x) = \frac{\hrev(\Delta;x)}{(1-x)^d},
\]
where $\hrev$ denotes the reverse $h$-polynomial:
\[
\hrev(\Delta;x) = x^dh\left(\Delta;\frac{1}{x}\right).
\]
\subsection{Shellable complexes}
Suppose $\Delta$ is a {\em pure} simplicial complex, meaning that all
facets (maximal simplices) have the same dimension $d-1$. A {\em
  shelling order} for $\Delta$ is a total ordering $F_1, \dots, F_t$
of the facets of $\Delta$ such that $F_j \cap
\left(\cup_{i<j}F_i\right)$ is pure of dimension $d-2$ for all $j=2, \dots, t$. We say that $\Delta$ is {\em shellable} if a shelling order for $\Delta$ exists.

One good reason to care about shellability is that it implies Cohen-Macaulayness.

\section{The objects of study} \label{se:defin}
Suppose $\h$ is a hyperplane arrangement in $\R^n$ such that $\cap \h
= \{0\}$. Then, $\h$ determines a regular cell decomposition $\dH$ of the unit sphere $S^{n-1}$. In short, each point $p$ on $S^{n-1}$ has an associated sign vector in $\{0,-,+\}^{|\h|}$ recording for each hyperplane $h\in \h$ whether $p$ is on, or on the negative, or on the positive side of $h$ (for some choice of orientations of the hyperplanes). A cell in $\dH$ consists of the set of points with a common sign vector. The face poset of $\dH$ is the big face lattice of the corresponding oriented matroid, see \cite{BLSWZ}. 

If $\dH$ is a simplicial complex, then $\h$ is called {\em simplicial}. A prime example of a simplicial hyperplane arrangement is the collection of reflecting hyperplanes of a finite Coxeter group. In this case, $\dH$ coincides with the Coxeter complex. 

From now on, let $\h$ be a simplicial hyperplane arrangement.

Consider an antichain $\A$ in $\LH$. We say that the subspace arrangement $\A$ is {\em embedded} in $\h$. Observe that $\cup \A \cap S^{n-1}$, which is known as the {\em link} of $\A$, has the structure of a simplicial subcomplex of $\dH$. This subcomplex is the principal object of study in this paper. We denote it $\dAH$.

\begin{example} \label{ex:graph}
{\em A graph $G = ([n],E)$ determines a {\em graphical hyperplane arrangement} $\G$ in the $(n-1)$-dimensional subspace of $\R^n$ given by the equation $x_1 + \dots + x_n = 0$. There is one hyperplane in $\G$ for each edge in $E$; the hyperplane corresponding to the edge $\{i,j\}$ has the equation $x_i = x_j$.

The arrangement $\Kn$ corresponding to the complete graph is nothing but the braid arrangement $\Sn$ which is simplicial. Any graph $G$ thus determines a simplicial complex $\dG$. It coincides with Steingr\'{\i}msson's coloring complex of $G$ which was denoted $\Delta_G$ in the Introduction. The complex $\dG$ also appeared under the name $\Delta_{m,J}$ in \cite{HRW}.}
\end{example}

We remark that the homotopy type of the link of $\A$, hence of $\dAH$,
can be computed in terms of the order complexes of lower intervals in
$L_\A$ by a formula of Ziegler and \v{Z}ivaljevi\'c \cite{ZZ}. When
$\A$ consists of hyperplanes we may simply note that $\dAH$ is
homotopy equivalent to the $(n-1)$-sphere with one point removed for
each connected region in the complement $\R^n\setminus \cup
\A$. Denoting by $R(\A)$ the number of such regions, $\dAH$ is thus
homotopy equivalent to a wedge of $R(\A)-1$ spheres of dimension $n-2$
in this case. For the arrangements $\G$ of Example \ref{ex:graph} it
is not difficult to see that $R(\G)$ equals the number $\text{AO}(G)$
of acyclic orientations of $G$. Thus, $\dG$ has the homotopy type of a
wedge of $AO(G)-1$ $(n-3)$-spheres (\cite{HRW,jonsson}). In
particular, the reduced Euler characteristic of $\dG$ is $\pm (\text{AO}(G)-1)$ (\cite[Theorem 17]{steingrimsson}).

\section{Shellability in the hyperplane case}\label{se:shellable}

Our goal in this section is to show that $\dAH$ is shellable whenever $\A$ consists of hyperplanes. Applied to the complexes $\dG$ of Example \ref{ex:graph} this answers affirmatively a question of Steingr\'{\i}msson \cite{steingrimsson} which was restated in \cite{jonsson}. The key tool is a particular class of shellings of $\dH$ determined by the {\em poset of regions} of $\h$ which we now define.

The complement $\R^n\setminus \cup \h$ is cut into disjoint open
regions by $\h$. Restricting to the unit sphere, their closures are the facets of $\dH$. Let $\f = \f(\h)$ be the set of such facets. For $R, R^\prime\in \f$, say that $h\in \h$ {\em separates} $R$ and $R^\prime$ if their respective interiors are on different sides of $h$. 

Choose a {\em base region} $B \in \f$ arbitrarily. We have a distance function $\ell:\f \to \N$ which maps a region $R$ to the number of hyperplanes in $\h$ which separate $R$ and $B$. Now, for two regions $R, R^\prime\in \f$, write $R \lhd R^\prime$ iff $R$ and $R^\prime$ are separated by exactly one hyperplane in $\h$ and $\ell(R) = \ell(R^\prime)-1$. The poset of regions $\ph$ is the partial order on $\f$ whose covering relation is $\lhd$. It was first studied by Edelman \cite{edelman}.

From the point of view of this paper, the most important property of
$\ph$ is the following.  
\begin{theorem}[Theorem 4.3.3 in \cite{BLSWZ}] \label{th:weak}
Any linear extension of $\ph$ is a shelling order for $\dH$.
\end{theorem}

We are now ready to state and prove the main result of this section.

\begin{theorem}\label{th:shelling}
If $\A$ consists of hyperplanes, then $\dAH$ is shellable.
\begin{proof}
We proceed by induction over $|\A|$. When $\A=\{A\}$, we may apply Theorem \ref{th:weak} since $\dAH = \Delta_{\h/A}$ in this case.

Now suppose $|\A|\geq 2$ and that we have a shelling order for $\Delta_{\A\setminus A,\h}$ for some $A\in \A$. We will append the remaining facets to this order. 

The remaining facets are the facets of $\Delta_{\{\A\},\h} = \Delta_{\h/A}$. They are divided into equivalence classes in the following way: $F$ and $G$ belong to the same class iff their interiors belong to the same connected component of $\R^n\setminus \cup(\A\setminus A)$ (or, equivalently, to the same connected component of $A\setminus \cup(\A/A)$). Observe that if $F$ and $G$ belong to different classes, then $F\cap G \in \Delta_{\A\setminus A,\h}$. Thus, it is enough to show that the facets in any equivalence class can be appended to the shelling order for $\Delta_{\A\setminus A,\h}$. 

Without loss of generality, consider the class which contains the maximal element in $P_{\h/A}$, i.e.\ the region opposite to the base region. Call this class $C$. If $F \in C$ and $G\not \in C$ for $F,G\in P_{\h/A}$, then some hyperplane in $\A/A \subseteq \h/A$ separates $F$ from $G$, and $G$ is on the positive side of this hyperplane. Thus, $F \not \leq G$. This shows that $C$ is an order filter in $P_{\h/A}$. According to Theorem \ref{th:weak}, $\Delta_{\h/A}$ has a shelling order which ends with the facets in $C$. Now observe that $(\cup C) \cap (\cup (P_{\h/A}\setminus C)) = (\cup C)\cap \Delta_{\A\setminus A,\h}$. The facets in $C$ may therefore be appended in this order to the shelling order for $\Delta_{\A\setminus A,\h}$.
\end{proof} 
\end{theorem}

\section{The $h$-polynomial of $\dAH$} \label{se:h-vector}

 For brevity we write $h(\A,\h;x)$ meaning $h(\dAH;x)$ and similarly for $\hrev$. The following result of Steingr\'{\i}msson serves as a motivating example for this section:
\begin{theorem}[Theorem 13 in \cite{steingrimsson}]\label{th:steingrimsson}
Recall the complex $\dG$ defined in Example \ref{ex:graph}. We have
\[
\frac{x\hrev(\G,\Sn;x)}{(1-x)^n} = \sum_{m\geq 0}\left(m^n-\chrom_G(m)\right)x^m,
\]
where $\chrom_G$ is the chromatic polynomial of $G$.
\end{theorem}

This theorem is interesting because of the connection between the left
hand side and the Hilbert series of the Stanley-Reisner ring
$k[\dG]$. In \cite{brenti}, Brenti began a systematic study of which
polynomials arise as Hilbert polynomials of standard graded
algebras. A question left open in \cite{brenti}, and later answered
affirmatively by Almkvist \cite{almkvist}, was whether chromatic
polynomials of graphs have this property. Theorem
\ref{th:steingrimsson} implies something similar, namely that
$(m+1)^n-\chrom_G(m+1)$ is the Hilbert polynomial (in $m$) of a standard graded algebra; for details, see Corollary \ref{co:Sn} below.

It is well-known that $\chrom_G(x) = x\charac(\G;x)$; one way to prove
it is to compare Theorem \ref{th:DR} with the standard
deletion-contraction recurrence for $\chrom_G$. The identity suggests the
possibility of extending Theorem \ref{th:steingrimsson} to other
complexes $\dAH$. This turns out to be possible at least if $\h \in
\{\Sn,\Bn\}$ and is the topic of this section. 

Given a subspace $T$ of $\R^n$, let $d(T)$ denote its dimension. For a subspace arrangement $\mathcal{T}$, we also write
\[
d(\mathcal{T}) = \max_{T\in \mathcal{T}}d(T).
\]

\begin{lemma}\label{le:recursion}
Let $A \in \A$. Then,
\[
\begin{split}
h(\A,\h;x) & = \,\,\,\,\,\, (x-1)^{d(\A) - d(\A\setminus A)}h(\A\setminus A,\h;x)\\
& \,\,\,\,\,\, +(x-1)^{d(\A)-d(A)}h(\{A\},\h;x)\\
& \,\,\,\,\,\, - (x-1)^{d(\A)-d(\A/A)}h(\A/A,\h/A;x).
\end{split}
\]
\begin{proof}
Each simplex in $\dAH$ belongs to $\Delta_{\A\setminus A,\h}$ or to $\Delta_{\{A\},\h}$ or to both. Also, observe that $\Delta_{\A\setminus A,\h} \cap \Delta_{\{A\},\h} = \Delta_{\A/A,\h/A}$. Denoting by $f_i(\mathcal{S},\mathcal{T})$ the number of $i$-dimensional simplices in $\Delta_{\mathcal{S},\mathcal{T}}$, we thus obtain for all $i$
\[
f_i(\A,\h) = f_i(\A\setminus A,\h) + f_i(\{A\},\h) - f_i(\A/A,\h/A).
\]
The lemma now follows from the fact that $\dim(\Delta_{\mathcal{S},\mathcal{T}}) = d(\mathcal{S})-1$.
\end{proof}
\end{lemma}
We may use Lemma \ref{le:recursion} to recursively compute $h(\A,\h;x)$. As it turns out, this recursion is particularly useful when $\h \in \{\Sn, \Bn\}$. The reason is given by the following two lemmata. 

\begin{lemma}\label{le:eulerian}
We have
\[
\frac{x\hrev(\Delta_{\Sn};x)}{(1-x)^{n+1}}=\sum_{m\geq 0}m^nx^m
\]
and
\[
\frac{\hrev(\Delta_{\Bn};x)}{(1-x)^{n+1}}=\sum_{m\geq 0}(2m+1)^nx^m.
\]
\begin{proof}
The complexes $\Delta_{\Sn}$ and $\Delta_{\Bn}$ coincide with the
Coxeter complexes of types $A_{n-1}$ and $B_n$, respectively. For the
$h$-polynomials this implies that $x\hrev(\Delta_{\Sn};x) = A_n(x)$
and $\hrev(\Delta_{\Bn};x)=B_n(x)$, where $A_n$ is the $n$th Eulerian
polynomial and $B_n$ is the $n$th $B$-Eulerian polynomial, see
\cite{brenti2}. The assertions are well-known properties of these
polynomials \cite[Theorem 3.4.ii]{brenti2}.  
\end{proof}
\end{lemma} 

\begin{lemma}\label{le:single}
$\,$
\begin{enumerate}
\item[(i)] For any subspace $A\in L_{\Sn}$, we have
\[
\frac{x\hrev(\{A\},\Sn;x)}{(1-x)^{d(A)+2}} = \sum_{m\geq 0}m^{d(A)+1}x^m.
\]
\item[(ii)] For any subspace $\A \in L_{\Bn}$, we have
\[
\frac{\hrev(\{A\},\Bn;x)}{(1-x)^{d(A)+1}} = \sum_{m\geq 0}(2m+1)^{d(A)}x^m.
\]
\end{enumerate}
\begin{proof}
A key property of $\Sn$ ($\Bn$), which is readily checked, is that its restriction to any subspace in the intersection lattice is again a type $A$ ($B$) hyperplane arrangement. Thus, $\Delta_{\{A\},\Sn} = \Delta_{\Sn/A} \cong \Delta_{\mathcal{S}_{d(A)+1}}$ ($\Delta_{\{A\},\Bn} = \Delta_{\Bn/A} \cong \Delta_{\mathcal{B}_{d(A)}}$). The assertions now follow from Lemma \ref{le:eulerian}
\end{proof}
\end{lemma}

The leading term of $\charac(\A;x)$ is always $x^n$, where $n$ is the dimension of the ambient space. It is convenient to introduce the {\em tail} $\tail(\A;x) = x^n-\charac(\A;x)$. 

When $\A$ consists of hyperplanes, the following result coincides with Theorem \ref{th:steingrimsson}.

\begin{theorem}\label{th:charac}\label{th:Sn}
Suppose $\A$ is a subspace arrangement embedded in $\Sn$. Then,
\[
\frac{x\hrev(\A,\Sn;x)}{(1-x)^{d(\A)+2}} = \sum_{m \geq 0}mT(\A;m)x^m.
\]
\begin{proof}
We proceed by induction over $|\A|$, noting that $|\A\setminus A|<|\A|$ and $|\A/A|<|\A|$ for every $A\in \A$. If $|\A| = 1$, we have $\charac(\A;m)=m^{n-1}-m^{d(\A)}$, so that $\tail(\A;m) = m^{d(\A)}$, and the theorem follows from part (i) of Lemma \ref{le:single}.

Now suppose $|\A|\geq 2$ and pick a subspace $A\in \A$. Using Lemma \ref{le:recursion} and the induction hypothesis, we obtain
\[
\begin{split}
\frac{x^{d(\A)+1}h\left(\A,\Sn;\frac{1}{x}\right)}{(1-x)^{d(\A)+2}} \, & =
\,\,\,\,\,\, \left(\frac{1-x}{x}\right)^{d(\A)-d(\A\setminus A)}\frac{x^{d(\A)+1}h\left(\A\setminus A,\Sn;\frac{1}{x}\right)}{(1-x)^{d(\A)+2}} \\
 &  \,\,\,\,\,\, + \left(\frac{1-x}{x}\right)^{d(\A)-d(A)}\frac{x^{d(\A)+1}h\left(\{A\},\Sn;\frac{1}{x}\right)}{(1-x)^{d(\A)+2}} \\
 &  \,\,\,\,\,\, - \left(\frac{1-x}{x}\right)^{d(\A)-d(\A/A)}\frac{x^{d(\A)+1}h\left(\A/A,\Sn/A;\frac{1}{x}\right)}{(1-x)^{d(\A)+2}} \\
 & = \,\,\,\,\,\, \sum_{m\geq 0}m(m^{n-1}-\charac(\A\setminus A;m))x^m \\
 & \,\,\,\,\,\, +  \sum_{m\geq 0}m(m^{n-1}-(m^{n-1}-m^{d(A)}))x^m\\
 & \,\,\,\,\,\, -  \sum_{m\geq 0}m(m^{d(A)}-\charac(\A/A;m))x^m\\
 & =  \,\,\,\,\,\, \sum_{m\geq 0}m(m^{n-1}-\charac(\A;m))x^m,
\end{split}
\]
where the last equality follows from Deletion-Restriction.

For completeness, we should also check the uninteresting case $|\A|=0$ which is not covered by the above arguments. Here, $\hrev(\emptyset,\Sn;x) = 0$ and $\tail(\emptyset;x)=0$, and the assertion holds. 
\end{proof}
\end{theorem}

Employing part (ii) of Lemma \ref{le:single} instead of part (i), and keeping track of the fact that $\Bn$ is an arrangement in $\R^n$, whereas $\Sn$ sits in $\R^{n-1}$, the proof of Theorem \ref{th:Sn} is easily adjusted to a proof of the next result.

\begin{theorem}\label{th:Bn}
Suppose $\A$ is a subspace arrangement embedded in $\Bn$. Then,
\[
\frac{\hrev(\A,\Bn;x)}{(1-x)^{d(\A)+1}} = \sum_{m \geq 0}\tail(\A;2m+1)x^m.
\]
\end{theorem} 

For subspace arrangements covered by Theorem \ref{th:Sn} or Theorem \ref{th:Bn}, we may now draw the promised algebraic conclusions. To this end, for a simplicial complex $\Gamma$ and a subcomplex $\Gamma^\prime \subseteq \Gamma$, let $\J_{\Gamma^\prime,\Gamma}$ be the ideal in the Stanley-Reisner ring $k[\Gamma]$ generated by the (equivalence classes of) monomials corresponding to simplices in $\Gamma$ that do not belong to $\Gamma^\prime$.

\begin{corollary}\label{co:Sn}
Suppose $\A$ is a subspace arrangement embedded in $\Sn$. Let $\Gamma$
denote the double cone over $\Delta_{\Sn}$, and write $\Gamma^\prime$
for the double cone over $\Delta_{\A,\Sn}$ with the same cone points. The following holds:
\begin{enumerate}
\item[(i)] The Hilbert polynomial of $k[\Gamma^\prime]$ is $F(k[\Gamma^\prime];m)=(m+1)\tail(\A;m+1)$.
\item[(ii)] The Hilbert polynomial of $\J_{\Gamma^\prime,\Gamma}$ is $F(\J_{\Gamma^\prime,\Gamma};m)=(m+1)\charac(\A;m+1)$.
\end{enumerate}
\begin{proof}
The dimension of $\Gamma^\prime$ is $d(\A)+1$. Taking a cone over a
simplicial complex does not affect the $\hrev$-polynomial. Thus,
\[
\hilb(k[\Gamma^\prime];x) = \frac{\hrev(\A,\Sn;x)}{(1-x)^{d(\A)+2}} = \frac{1}{x}\sum_{m\geq 0}m\tail(\A;m)x^m,
\] 
where the second equality follows from Theorem \ref{th:Sn}. This proves (i).

For (ii), we use that 
\[
k[\Gamma^\prime] \cong k[\Gamma]/\J_{\Gamma^\prime,\Gamma}.
\]
For the Hilbert series, this implies
\[
\hilb(k[\Gamma^\prime];x) = \hilb(k[\Gamma];x) - \hilb(\J_{\Gamma^\prime,\Gamma};x).
\]
From part (i) and the fact that
\[
\hilb(k[\Gamma])=\frac{\hrev(\Delta_{\Sn};x)}{(1-x)^{n+1}} = \frac{1}{x}\sum_{m\geq 0}m^nx^m, 
\]
we conclude
\[
\hilb(\J_{\Gamma^\prime,\Gamma};x) = \frac{1}{x}\sum_{m\geq 0}m^nx^m - \frac{1}{x}\sum_{m\geq 0}m\tail(\A;m)x^m = \frac{1}{x}\sum_{m\geq 0}m\charac(\A;m)x^m.
\]
\end{proof}
\end{corollary}
The situation for $\Bn$ is analogous, although we use cones instead of double cones. This is a manifestation of the fact that $\Bn$ and $\Sn$ differ by one in dimension.  
\begin{corollary}\label{co:Bn}
Suppose $\A$ is a subspace arrangement embedded in $\Bn$. Let $\Gamma$
denote the cone over $\Delta_{\Bn}$, and write $\Gamma^\prime$ for the
cone over $\Delta_{\A,\Bn}$ with the same cone point. Then, the following holds:
\begin{enumerate}
\item[(i)] The Hilbert polynomial of $k[\Gamma^\prime]$ is $F(k[\Gamma^\prime];m)=\tail(\A;2m+1)$.
\item[(ii)] The Hilbert polynomial of $\J_{\Gamma^\prime,\Gamma}$ is $F(\J_{\Gamma^\prime,\Gamma};m)=\charac(\A;2m+1)$.
\end{enumerate}
\begin{proof}
Proceeding as in the proof of Corollary \ref{co:Sn}, using Theorem
\ref{th:Bn} instead of Theorem \ref{th:Sn}, we prove (i) by observing
\[
\hilb(k[\Gamma^\prime];x) = \frac{\hrev(\A,\Bn;x)}{(1-x)^{d(\A)+1}} = \sum_{m\geq 0}\tail(\A;2m+1)x^m.
\]

For (ii), note that
\[
\hilb(k[\Gamma];x) = \frac{\hrev(\Delta_{\Bn};x)}{(1-x)^{n+1}} = \sum_{m\geq 0}(2m+1)^nx^m.
\]
Thus, 
\[
\hilb(\J_{\Gamma^\prime,\Gamma};x) = \sum_{m\geq 0}(2m+1)^nx^m - \sum_{m\geq 0}\tail(\A;2m+1)x^m = \sum_{m\geq 0}\charac(\A;2m+1)x^m.
\]
\end{proof}
\end{corollary}

Any hypergraph (without inclusions among edges) $G$ on $n$ vertices
corresponds to a subspace arrangement $\G$ embeddable in
$\Sn$. The construction is virtually the same as in Example \ref{ex:graph}; with the
hyperedge $\{i_1, \dots, i_t\}$ is associated the subspace given by
$x_{i_1}=\dots = x_{i_t}$. As for ordinary graphs (the hyperplane case), we
have $x\charac(\G;x)=\chrom_G(x)$, cf.\ \cite[Theorem 3.4]{stanley}. In
this way, Corollary 
\ref{co:Sn} allows us to interpret chromatic polynomials of
hypergraphs in terms of Hilbert polynomials. For ordinary graphs, this is the content of Steingr\'{\i}msson's \cite[Corollary 10]{steingrimsson}.

Corollary \ref{co:Bn}, too, has an impact on chromatic polynomials. Any {\em signed graph} (in the sense of Zaslavsky \cite{zaslavsky}) $G$ on $n$ vertices corresponds to a hyperplane arrangement $\G \subseteq \Bn$, and vice versa. A signed graph $G$ has a chromatic polynomial $\chrom_G(x)$, and $\chrom_G(x) = \charac(\G;x)$ \cite{zaslavsky2}.


\begin{thebibliography}{AB}
\bibitem{almkvist} G.\ Almkvist, The chromatic polynomial is a Hilbert polynomial, preprint 1998.

\bibitem{BLSWZ} A.\ Bj\"orner, M.\ Las Vergnas, B.\ Sturmfels, N.\ White and G.\ M.\ Ziegler, {\em Oriented matroids. 2nd ed.}, Encyclopedia of Mathematics and its Applications {\bf 46}, Cambridge Univ.\ Press, Cambridge, 1999. 

\bibitem{brenti} F.\ Brenti, Hilbert polynomials in combinatorics, {\em J.\ Algebraic Combin.\ }{\bf 7} (1998), 127--156.

\bibitem{brenti2} F.\ Brenti, $q$-Eulerian polynomials arising from Coxeter groups, {\em European J.\ Combin.\ }{\bf 15} (1994), 417--441.

\bibitem{edelman} P.\ H.\  Edelman, A partial order on the regions of $R^n$ dissected by hyperplanes, {\em Trans.\ Amer.\ Math.\ Soc.\ }{\bf 283} (1984), 617--631.

\bibitem{HRW} J.\ Herzog, V.\ Reiner and V.\ Welker, The Koszul property in affine semigroup rings, {\em Pacific J.\ Math.\ }{\bf 186} (1998), 39--65.

\bibitem{jonsson} J.\ Jonsson, The topology of the coloring complex,
  {\em J.\ Algebraic Combin.\ }{\bf 21} (2005), 311--329.

\bibitem{rota} G.-C.\ Rota, On the foundations of combinatorial
  theory. I. Theory of M\"obius functions, {\em Z.\
  Wahrscheinlichkeitstheorie und Verw.\ Gebiete} {\bf 2} (1964), 340--368.

\bibitem{stanley} R.\ P.\ Stanley, Graph colorings and related
  symmetric functions: ideas and applications: a description of
  results, interesting applications, \& notable open problems, {\em
  Discrete Math.\ }{\bf 193} (1998), 267--286.

\bibitem{steingrimsson} E.\ Steingr\'{\i}msson, The coloring ideal and coloring complex of a graph, {\em J.\ Algebraic Combin.\ }{\bf 14} (2001), 73--84.

\bibitem{zaslavsky} T.\ Zaslavsky, The geometry of root systems and signed graphs, {\em Amer.\ Math.\ Monthly} {\bf 88} (1981), 88--105.

\bibitem{zaslavsky2} T.\ Zaslavsky, Signed graph coloring, {\em Discrete Math.\ }{\bf 39} (1982), 215--228.

\bibitem{ZZ} G.\ M.\ Ziegler and R.\ T.\ \v{Z}ivaljevi\'c, Homotopy types of subspace arrangements via diagrams of spaces, {\em Math.\ Ann.\ }{\bf 295} (1993), 527--548. 
\end{thebibliography}
\end{document}